\documentclass [10pt,oneside,a4paper,mathscr]{amsart}

\usepackage{amscd,amsmath,amssymb,euscript}
\usepackage[frame,cmtip,curve,arrow,matrix,line,graph]{xy}

\date{}

\title{Quantum groups via Hall algebras of complexes}
\author{Tom Bridgeland}
\address{All Souls College, Oxford OX1 4AL, U.K.}
\renewcommand{\footnote}[1]{}
\newtheorem{thm}{Theorem}[section]

\newtheorem{lemma}[thm]{Lemma}
\newenvironment{pf}{\paragraph{Proof}}{\qed\par\medskip}
\theoremstyle{definition}

\newtheorem{thm*}[thm]{Theorem$^*$}

\newtheorem*{example*}{Example}

\newcommand {\A}{\mathcal A}

\newcommand{\D}{\operatorname{D}}
\newcommand{\diam}{\diamond}
\newcommand {\C}{\mathbb C}

\newcommand {\Hom}{\operatorname{Hom}}
\newcommand {\Aut}{\operatorname{Aut}}
\newcommand {\Ext}{\operatorname{Ext}}

\renewcommand {\H}{\operatorname{\mathcal H}}
\newcommand{\Vect}{\operatorname{Vect}}
\newcommand {\id}{\operatorname{id}}
\newcommand{\cplx}[4] {\begin{xy}
\xymatrix@=0pt@C=20pt@R=15pt{#1
\ar@<.5ex>[r]^{#2}		& \ar@<.5ex>[l]^{#3}  #4}
\end{xy}  }
\newcommand{\mat}[4]{\begin{pmatrix}#1&#2\\#3&#4\end{pmatrix}}

\newcommand{\cl}[1]{{\hat{#1}}}

\newcommand{\gen}[1]{[#1]}
\newcommand{\R}{{\operatorname{R}}}
\newcommand{\lra}{\longrightarrow}

\newcommand {\<}{\langle}
\renewcommand {\>}{\rangle}
\newcommand {\Omit}[1]{}

\newcommand{\isom}{\cong}

\newcommand{\tensor}{\otimes}
\renewcommand{\dag}{*}
\newcommand{\fb}{\mathfrak{b}}
\newcommand{\e}{\operatorname{e}}

\newcommand{\Zt}{{\mathbb{Z}_2}}
\newcommand{\Ho}{\operatorname{Ho}}
\renewcommand{\P}{\mathcal{P}}
\newcommand{\CC}{\mathcal{C}}
\newcommand{\Z}{\mathbb{Z}}
\newcommand{\red}{\operatorname{red}}
\newcommand{\g}{\mathfrak{g}}
\newcommand{\n}{\mathfrak{n}}
\newcommand{\h}{\mathfrak{h}}
\newcommand{\F}{\mathbb{F}}
\newcommand{\blank}{\kern.0mm\smash{-}\kern.0mm}
\newcommand{\blob}{{\scriptscriptstyle\bullet}}

\newcommand{\E}{E}
\newcommand{\K}{K}

\newcommand{\longinto}{\lhook\joinrel\longrightarrow}

\newcommand{\lRa}[1]{\stackrel{#1}{\lra}}

\newcommand{\tw}{\operatorname{tw}}
\newcommand{\Htw}{\H_{\tw}}
\renewcommand{\DH}{\operatorname{\mathcal{DH}}}
\newcommand{\DHred}{\operatorname{\mathcal{DH}_{red}}}
\newcommand{\U}{\operatorname{U}}
\newcommand{\im}{\operatorname{im}}

\newcommand{\Iso}{\operatorname{Iso}}

\begin{document}
\begin{abstract}{We describe quantum enveloping algebras of symmetric Kac-Moody Lie algebras via a finite field Hall algebra construction involving $\Z_2$-graded complexes of quiver representations.}\end{abstract}
\maketitle
\section{Introduction}

\subsection{Background}Let $Q$ be a quiver without oriented cycles, and let $\A$ be the category of finite-dimensional representations of $Q$ over the finite field $k=\F_q$. The graph underlying $Q$ determines a generalized Cartan matrix, and we let $\g=\n^+\oplus\h\oplus \n^-$ denote the corresponding  derived Kac-Moody Lie algebra.
A famous theorem of Ringel \cite{ringel}, as extended by Green \cite{green},   gives an embedding of associative algebras
\[R\colon \U_{t}(\n^+) \longinto\Htw(\A).\]
On the left is the positive part of the quantum  enveloping algebra of $\g$, specialized at $t=\sqrt{q}$.  On the right is  the twisted Hall algebra  of the abelian category $\A$.
The map $R$ is an isomorphism precisely when the graph underlying $Q$ is  a simply-laced Dynkin diagram.

In the general case,  a generator $E_i$ of the quantum group is mapped to the generator of $\Htw(\A)$ defined by the simple module at the corresponding vertex of
$Q$. Thus the algebra $\U_{t}(\n^+)$  can  be recovered as the subalgebra of $\Htw(\A)$ generated by these elements.  Green later showed  how to give a categorical construction of the quantum coproduct \cite{green}, and Xiao  described the antipode \cite{xiao}. Thus the complete Hopf algebra structure of $\U_t(\n^+)$ can be described in terms of the  abelian category $\A$. We recommend \cite{schiffmann} for a readable survey of these results.

There have been various attempts to extend this picture so as  to give a Hall algebra description of the full quantum group $\U_t(\g)$ (see for example \cite{kapranov}).
It was realised early on that $\U_t(\g)$  should correspond in some way to $\Zt$-graded complexes in $\A$, with the pieces $\U_t(\n^+)$ and $\U_t(\n^-)$  of the triangular decomposition corresponding to complexes concentrated in degree 0 and 1 respectively. It was also understood that $\U_t(\h)$  corresponds naturally to the group algebra of the Grothendieck group $K(\A)$. Without a  unified treatment of these different  pieces  however, it is difficult  to see how to get  relations of the form
\begin{equation}\label{den}[E_i,F_i]=\frac{K_i-K_i^{-1}}{t-t^{-1}}\end{equation}
which hold in $\U_t(\g)$.

Much of the work on this problem has focused on the derived category of $\A$. There is then a problem, because  the non-uniqueness of cones in a triangulated category prevents the obvious definition of the Hall algebra from being associative. To{\"e}n \cite{toen} solved this problem using a dg-enhancement, and it was later proved  that the resulting formula defines an associative multiplication using only the axioms of a triangulated category \cite{xiao_xu}.  But in both cases, the finiteness assumptions required rule out the case of derived categories of $\Z_2$-graded complexes.
Peng and Xiao \cite{peng_xiao,peng_xiao2} did manage to use $\Z_2$-graded complexes in $\A$ to construct the full  Lie algebra $\g$, though the Cartan elements appear in this construction in a rather ad hoc manner, and the crucial relations
$[E_i,F_i]=H_i$
have to be put in by hand.

In this paper we show how to construct the full quantum group $U_t(\g)$ from the Hall algebra of the abelian category of $\Z_2$-graded complexes in $\A$. We also establish the connection between this Hall algebra of complexes $\DH(\A)$ and the Drinfeld double of the extended Hall algebra of $\A$.

\subsection{Main result}

Let $\A$ be an abelian category over a finite field $k=\mathbb{F}_q$, with finite-dimensional morphism spaces, and let $\P\subset\A$ be the subcategory of projective objects.
Let  $\CC(\A)$ be the abelian category of  $\Z_2$-graded complexes in $\A$
\[\cplx{M_1}{f}{g}{M_0}, \qquad f\circ g = 0 = g\circ f,\] 
and let $\CC(\P)$ be the subcategory of projective complexes.   We consider the Hall algebra of the category $\CC(\P)$, and twist the multiplication using the Euler form of $\A$. Our algebra $\DH(\A)$ is then the localization of this twisted Hall algebra with respect to the set of acyclic complexes
\[\DH(\A)=\Htw(\CC(\P))\big[ \gen{M_\blob}^{-1}: H_*(M_\blob)=0\big].\]
We also consider   a reduced version $\DHred(\A)$. To define this, first note that the shift  functor defines an involution of $\CC(\P)$;  we write it as $M_\blob \leftrightarrow M_\blob^\dag$. Now define
\[\DH_{\red}(\A)=\DH(\A)\big/(\gen{M_\blob}-1: H_*(M_\blob)=0,\  M_\blob\isom M_\blob^\dag),\]
or in other words, impose relations  $\gen{M_\blob}=1$ for each acyclic complex $M_\blob$ that is invariant under shift.

%

\begin{thm}
\label{first}
Let $\A$ be the category of finite-dimensional $\F_q$-representations of a finite quiver $Q$ without oriented cycles. Then there is an embedding of algebras
\[R\colon \U_{t}(\g) \longinto\DH_{\red}(\A),\]
where $t=\sqrt q,$ and $\g$ is the derived\footnote{Understand this derived business} Kac-Moody Lie algebra associated to $Q$. The map $R$ is an isomorphism precisely when the underlying graph of $Q$ is a simply-laced Dynkin diagram.
\end{thm}

\smallskip

In the case when $Q$ is the  $A_1$ quiver, and $\A$ is the category of finite-dimensional vector spaces over $k$, we get an isomorphism
\[R\colon \U_{t}(\mathfrak{sl_2})\lra\DH_{\red}(\A).\]
Up to scalar factors, the generators correspond to complexes as follows
\[E_i : \cplx{0}{0}{0}{k}, \quad F_i : \cplx{k}{0}{0}{0}\]
\[K_i:\cplx{k}{1}{0}{k}, \quad K_i^{-1} : \cplx{k}{0}{1}{k}.\]
Note that the elements $K_i$ and $K_i^{-1}$ correspond to acyclic complexes, which become zero objects on passing to the derived category.
This example was known to Kapranov as long ago as 1995.

\subsection{Assumptions}
\label{ass}
Throughout the paper $\A$ denotes an abelian category. Let us now discuss the various conditions we impose on this category. 
We shall always assume that 
\smallskip\begin{itemize}
\item[(a)] $\A$ is essentially small, with finite morphism spaces,
\item[(b)] $\A$ is linear over $k=\F_q$,
\item[(c)] $\A$ is of finite global dimension and has enough projectives.
\end{itemize}

The essentially small condition is the statement that $\A$ is equivalent to a category whose objects form a set. Assumption (a) seems to be essential to the naive approach to Hall algebras involving counting isomorphism classes. It also ensures that $\A$ is a Krull-Schmidt category.
In contrast, assumption  (b) is not really necessary anywhere in  the paper, and is made for aesthetic reasons. Without it one has to use multiplicative Euler forms instead of additive ones, as for example in \cite{schiffmann}. 

Assumption (c) is equivalent to the statement that every object has a finite projective resolution. For the material of Section 3, this is not really needed, and if one replaces  $K(\A)$ by the Grothendieck group of the category  of projectives $\P\subset \A$, the definitions can all be made in the same way.  It is however questionable whether the resulting Hall algebra is the correct object to consider in these more general cases.

From Section 4 onwards we shall also assume that
\smallskip
\begin{itemize}
\item[(d)] $\A$ is hereditary, that is of global dimension at most 1,
\item[(e)] nonzero objects in $\A$ define nonzero classes in $K(\A)$.
\end{itemize}

To explain condition (e) first define $K_{\geq 0}(\A)\subset K(\A)$ to be the positive cone in the  Grothendieck group, consisting of classes of objects of $\A$, rather than formal differences of such. Then (e) is  equivalent to the statement that the rule
\begin{equation}
\label{order}\alpha\leq \beta \iff \beta-\alpha\in K_{\geq 0}(\A)\end{equation}
defines a partial order on $K(\A)$.  It holds, for example, if $\A$ has finite length. 

The main  example we have in mind, which satisfies all these conditions,  is the category of finite-dimensional representations of a finite quiver with no oriented cycles.

\subsection{Further directions}

In the hereditary case the Drinfeld double construction has been frequently used as a substitute for the missing derived Hall algebra \cite{burban, burb_schiff}. Xiao \cite{xiao} first used this construction as an indirect method to 
construct the full quantum enveloping algebra $\U_t(\g)$  from the abelian category of representations of $Q$.  Given an abelian category satisfying conditions (a)--(e), work of Ringel \cite{ringel} and Green \cite{green}  defines a self-dual topological bialgebra $\H_{\tw}^{\e}(\A)$ known as the extended twisted Hall algebra. It is possible to prove
\begin{thm}
\label{thm}
Suppose $\A$ is an abelian category satisfying the conditions (a)--(e) and which is either artinian or noetherian.  Then the algebra $\DH(\A)$ is isomorphic (as an associative algebra)  to the Drinfeld double of the bialgebra $\H_{\tw}^{\e}(\A)$.
Similarly, the  algebra  $\DH_{\red}(\A)$ is isomorphic to the reduced Drinfeld double.\end{thm}

Details, as they say, will appear elsewhere. It follows from Theorem \ref{thm} that, in the hereditary case, under suitable finiteness assumptions, the Hall algebra $\DH(\A)$ has the structure of a bialgebra.
 We leave the detailed categorical description of this structure for future research.
We also leave for future research the question of derived invariance. Recent results of Cramer \cite{cramer} together with Theorem \ref{thm} suggest that for hereditary categories the algebra $\DH(\A)$ should be functorial with respect to derived equivalences. I have not yet managed to give a direct construction of  this action within the framework of this paper. 

There are several other lines of research suggested by the results of this paper. One obvious thought is to generalize away from Hall algebras defined over finite fields. In particular, it should be possible to consider Hall algebras of perverse sheaves in the style of Lusztig. This has the potential to lead to new results on canonical bases. One could also consider Hall algebras  of constructible functions in characteristic zero, which should give rise  to (generalized)  Kac-Moody Lie algebras. This approach should be closely related to work of Peng and Xiao \cite{peng_xiao,peng_xiao2}.

Perhaps the most interesting generalization would be to consider categories of global dimension $\geq 2$. The definition of $\DH(\A)$ given in this paper makes perfect sense in this case, but I have very little idea about the properties of the resulting algebras, since the crucial results of Section \ref{her}  rely heavily on the hereditary assumption. A good test case would be to consider the category of finitely-generated modules over the polynomial algebra $\C[x,y]$. One might hope that this would give a link with the results of Grojnowski and Nakajima.

\subsection*{Acknowledgements}
What little I know about quantum groups I learned from  Valerio Toledano Laredo, and discussions with him played a crucial role in the early stages of this project. I'm very grateful to him for his many careful explanations over the past few years. I'd also like to thank Igor Burban for his helpful comments on the manuscript, and particularly for pointing me to  \cite[Remark 3.12]{burban}, which simplified the exposition considerably.

\subsection*{Notation}
 We fix a  field $k=\F_q$ with $q$ elements, and set $t=+\sqrt q$.  We write $\Iso(\A)$ for the set of isomorphism classes of a small category $\A$. The symbol $|S|$ denotes the number of elements of a finite set $S$.

\section{Background}
In this section we recall the basic definitions concerning quantum enveloping algebras and  Hall algebras,  and state Ringel's theorem. Almost all of this material can be found in Schiffmann's notes \cite{schiffmann}.
\subsection{Quantum enveloping algebra} 

Let $\Gamma$ be a finite graph with vertex set $\{1,\cdots,n\}$. Let $n_{ij}$ be the number of edges connecting vertex $i$ and $j$. We always assume that $\Gamma$ has no loops, or equivalently that $n_{ii}=0$ for all $i$. Let
\[a_{ij}=2\delta_{ij}-n_{ij}.\]
Then $a_{ij}$ is a symmetric generalised Cartan matrix, which we refer to as the Cartan matrix of the graph $\Gamma$. 
Let $\U_t(\g)$ denote the  quantum enveloping algebra of the corresponding derived Kac-Moody Lie algebra, specialised at $t$. 
As an associative algebra $\U_t(\g)$ is generated by symbols
$E_i, F_i, K_i, K_i^{-1}$ 
subject to the  relations\smallskip
\begin{equation} \label{firs}K_i* K_i^{-1}=1, \quad K_i^{-1}* K_i=1, \qquad [K_i, K_j] = 0, \end{equation}
\begin{equation} K_i * E_j = t^{ a_{ij} }\cdot  E_j * K_i, \qquad K_i * F_j = t^{- a_{ij} }\cdot  F_j * K_i,\end{equation}
\begin{equation}\label{anom} [E_i,F_i] =  {\frac{K_i-K_i^{-1}}{t-t^{-1}}}, \qquad [E_i,F_j]=0 \text { for }i\neq j,\end{equation}
and also the quantum Serre relations
\begin{equation}\label{serre1}\big.\sum_{n=0}^{1-a_{ij}} (-1)^n {{1-a_{ij}}\brack{n}}\cdot E_i^n * E_j * E_i^{1-a_{ij}-n} = 0, \quad i\neq j,\end{equation}
\begin{equation}\label{serre2}\big.\sum_{n=0}^{1-a_{ij}} (-1)^n {{1-a_{ij}}\brack{n}} \cdot  F_i^n * F_j * F_i^{1-a_{ij}-n} = 0,\quad i\neq j\end{equation}
where the coefficients
\begin{equation} {n\brack m}={ \frac{[n]! \cdot [m]!} {[n+m]!}}, \quad [n]!=\prod_{i=1}^n [i], \quad [n]=\frac{t^n-t^{-n}}{t-t^{-1}},\end{equation}
are quantum binomials. Note that there is an involution
\[U_t(\g) \stackrel{\sigma}{\longleftarrow\joinrel\longrightarrow}U_t(\g),\]
 defined on generators by $ E_i\stackrel{\sigma}{\longleftrightarrow} F_i$ and $K_i\stackrel{\sigma}{\longleftrightarrow} K_i^{-1}$.

\subsection{Triangular decomposition}

Define subalgebras
\[\U_t(\n^+), \U_t(\h), \U_t(\n^-)\subset \U_t(\g)\]
generated by the $E_i$, the $K_i^{\pm 1}$, and the $F_i$ respectively. Then
\[\U_t(\h)\isom \C[x_1^{\pm1},\cdots, x_n^{\pm 1}],\]
and the  involution $\sigma$ identifies $\U_t(\n^+)$ and $\U_t(\n^-)$.
The following result is referred to as the triangular decomposition of $\U_t(\g)$.

\begin{lemma}\label{tri}The multiplication map
\[\U_t(\\n^+)\tensor_\C \U_t(\h)\tensor_\C \U_t(\n^-)\lra \U_v(\g)\]
 is 
an isomorphism.
\end{lemma}

\begin{pf}This is standard. See for example \cite[Cor. 3.2.5]{lusztig}.
\end{pf}

We also define \[\U_t(\fb^+),\U_t(\fb_-)\subset \U_t(\g)\]
 to be the subalgebras generated by $E_i,K_i^{\pm 1}$ and by $F_i,K_i^{\pm 1}$, respectively. The multiplication maps 
\[\U_t(\h)\tensor_\C \U_t(\n^\pm)\lra \U_t(\fb^{\pm})\]
are isomorphisms of vector spaces.

\subsection{Hall algebra}
\label{hall}Let  $\A$ be an abelian category satsifying the assumptions (a)--(c) of Section \ref{ass}. Given objects $A,B,C\in \A$, define \[\Ext^1_\A(A,C)_B\subset \Ext^1_\A(A,C)\]
 to be the subset parameterising 
extensions with middle term isomorphic to $B$.
We define the Hall algebra $\H(\A)$ to be the vector space with basis  indexed by elements $A\in\Iso(\A)$, and with associative multiplication
defined by
\begin{equation}
\label{defprod}
\gen{A}\diam \gen{C}= \sum_{B\in\Iso(\A)} \frac{|\Ext^1_\A(A,C)_{B}|}{|\Hom_\A(A,C)|} \cdot \gen{B}.\end{equation}
The unit is $\gen{0}$.
Another form of the Hall product is perhaps more familiar. Given objects $A,B,C\in \A$, define numbers
\[g_{A,C}^B=\big|\big\{ C'\subset B :  C'\isom C,\ B/C'\isom A\big\}\big|.\]
Then one can prove\footnote{Look up reference} that
\[g_{A,C}^B=\frac{|\Ext^1_\A(A,C)_{B}|}{|\Hom_\A(A,C)|}\cdot\frac{|\Aut_\A(B)|}{|\Aut_\A(A)|\cdot|\Aut_\A(C)|}.\]
Thus in terms of alternative generators \[[[A]]=\frac{\gen{A}}{|\Aut_\A(A)|}\] the product takes the form
\[[[A]]\diam [[C]]= \sum_{B\in\Iso(\A)} g_{A,C}^B\cdot [[B]],\]
which is the definition used in \cite{ringel} for example.
 
\subsection{Euler form}
Let $K(\A)$ denote the Grothendieck group of $\A$. We write $\cl{A}\in K(\A)$ for the class of an object $A\in \A$. Let  $\K_{\geq 0}(\A)\subset K(\A)$ be the subset consisting of these classes  (rather than formal differences of them).
For objects $A,B\in \A,$ define
\[\<A,B\>=\big.\sum_{i\in\Z} (-1)^i \dim_k\Ext^i_\A(A,B).\]
The sum is finite by our assumptions on $\A$, and descends to give a bilinear form
\[\<\blank,\blank\>\colon K(\A)\times K(\A)\lra \Z\]
known as the Euler form.
Note that
\[\<P,Q\> =  \dim_k \Hom_\A(P,Q)\]
whenever $P\in \P$. 
We also consider the symmetrised form
\[(\blank,\blank)\colon K(\A)\times K(\A)\lra \Z,\]
defined by $( \alpha,\beta) =\<\alpha,\beta\> +\<\beta,\alpha\>$.

\subsection{Twisted and extended versions}
The twisted Hall algebra $\Htw(\A)$ is the same vector space as $\H(\A)$  but with twisted multiplication defined by
\[\gen{A}*\gen{C}=t^{\< \cl{A},\cl{C}\rangle}\cdot \gen{A}\diam \gen{C}.\]
Note that $\Htw(\A)$  has a grading
\[\Htw(\A)=\bigoplus_{\alpha\in K_{\geq 0}(\A)} \H_\alpha(\A),\]
where $\H_\alpha(\A)$ is the subspace spanned by elements $\gen{A}$ with $\cl{A}=\alpha$.
 We can  form the extended Hall algebra $\H_{\tw}^{\e}(\A)$ by adjoining symbols $K_\alpha$ for classes $\alpha\in K(\A)$, and imposing relations
\begin{equation}
\label{twisted}K_\alpha * K_\beta=K_{\alpha+\beta}, \quad K_\alpha * \gen{B} = t^{( \alpha,\cl{B})} \cdot \gen{B} * K_\alpha.\end{equation}
Note that  $\H_{\tw}^{\e}(\A)$ has a vector space basis consisting of the elements $K_\alpha * \gen{B}$ for  $\alpha\in K(\A)$ and $B\in \Iso(\A)$. 

\subsection{Ringel's Theorem}
Let $Q$ be a finite quiver with vertices $\{1,\cdots,n\}$ and no loops. Let $U_t(\g)$ be the quantized enveloping algebra corresponding to the Cartan matrix of the graph underlying $Q$.   Let $\A$ be the abelian category of  finite-dimensional representations of $Q$ over  $\F_q$, and for each vertex  $i$, let $S_i\in \A$ be the corresponding one-dimensional simple representation.

\begin{thm}[Ringel, Green]\label{ring}  There are injective homomorphisms of algebras
\[R\colon \U_t(\n^+)\longinto \Htw(\A),\quad R^{\e}\colon \U_t(\fb^+)\longinto \H_{\tw}^{\e}(\A),\]
defined on generators by
\[R(E_i)=R^{\e}(E_i)=\gen{S_i}/(q-1),\quad R^{\e}(K_i)=K_{\cl{S_i}}.\]
These maps are isomorphisms precisely if the underlying graph of $Q$ is  a simply-laced Dynkin diagram.
\end{thm}

In fact Ringel \cite{ringel} considered only the case of a Dynkin quiver and proved that $R$ is then an isomorphism. Green \cite{green} observed that Ringel's proof gave the existence of the algebra homomorphism $R$ in general, and used results of Lusztig to prove that it is an embedding.
For an account of the proof of Theorem \ref{ring} see \cite[Section 3.3]{schiffmann}.


%

\section{Hall algebras of complexes}
\label{reduction}
In this section we introduce Hall algebras of $\Zt$-graded complexes. Throughout $\A$ is an abelian category satisfying the assumptions (a)--(c) of Section \ref{ass}. 

\subsection{Categories of complexes}
\label{acy}

Let $\CC_{\Zt}(\A)$ be the abelian category of $\Zt$-graded complexes in $\A$.  The objects of this category  consist of  diagrams
 \begin{equation}
\label{complex}
\cplx{M_1}{d_1}{d_0}{M_0}\end{equation} in $\A$ such that $d_{i+1}\circ d_i=0$. A morphism $s_\blob\colon M_\blob\to \tilde{M}_\blob$ consists of a diagram
\[\begin{xy}
\xymatrix{
M_1\ar[d]_{s_1}\ar@<.5ex>[r]^{d_1}		& \ar@<.5ex>[l]^{d_0} M_0\ar[d]^{s_0}\\
\tilde{M}_1\ar@<.5ex>[r]^{\tilde{d}_1}	& \ar@<.5ex>[l]^{\tilde{d}_0}	 \tilde{M}_0}
\end{xy} \]
with $s_{i+1}\circ d_{i} = \tilde{d}_i\circ s_{i}$.
Two morphisms $s_\blob, t_\blob\colon M_\blob\to \tilde{M}_\blob$ are homotopic  if there are morphisms $h_i\colon M_i\to \tilde{M}_{i+1}$
such that
\[t_i-s_i=\tilde{d}_{i+1}\circ h_{i}+h_{i+1}\circ d_{i}.\]
We write $\Ho_\Zt(\A)$ for the category obtained from $\CC_\Zt(\A)$ by identifying homotopic morphisms.
The  shift  functor defines an involution \[\CC_\Zt(\A)\stackrel{\dag}{\longleftrightarrow}\CC_\Zt(\A)\]  which shifts the grading and changes the sign of the differential
\[ \begin{xy}
\xymatrix@=0pt@C=20pt@R=15pt{
M_1\ar@<.5ex>[r]^{d_1}		& \ar@<.5ex>[l]^{d_0}  M_0}
\end{xy}  \stackrel{\dag}{\longleftrightarrow}  \begin{xy}
\xymatrix@=0pt@C=20pt@R=15pt{
M_0\ar@<.5ex>[r]^{-d_0}		& \ar@<.5ex>[l]^{-d_1}  M_1}
\end{xy}\]
Note that every complex $M_\blob\in \CC_\Zt(\A)$ defines a class
\[\cl{M_\blob}=\cl{M_0}-\cl{M_1}\in K(\A).\]
We shall be mostly interested in the full subcategories \[\CC_\Zt(\P)\subset\CC_\Zt(\A),\quad\Ho_\Zt(\P)\subset \Ho_\Zt(\A),\] whose objects are complexes of projectives in $\A$.
After the next subsection we will drop the $\Zt$ subscripts and just write $\CC(\A)$ and $\Ho(\A)$.

\subsection{Root category}
Let $\D^b(\A)$ denote the ($\Z$-graded) bounded derived category of $\A$, with its shift functor $[1]$. Let $\R(\A)=\D^b(\A)/[2]$ be the orbit category, also known as the root category of $\A$. This has the same objects as $\D^b(\A)$, but the morphisms are given by
\[\Hom_{\R(\A)} (X,Y)=\bigoplus_{i\in \Z} \Hom_{\D^b(\A)}(X,Y[2i]).\]
 Since $\A$ is assumed to be of finite global dimension with  enough projectives  the category $\D^b(\A)$ is equivalent to the ($\Z$-graded) bounded homotopy category $\Ho^b(\P)$. Thus we can equally well define $\R(\A)$ as the orbit category of $\Ho^b(\P)$.
 
 \begin{lemma}
 \label{missing}
 There is a fully faithful functor
\[D\colon \R(\A)\lra \Ho_{\Zt}(\P)\]
sending a $\Z$-graded complex of projectives $(P_i)_{i\in\Z}$ to the $\Zt$-graded complex
\[\cplx{\bigoplus_{i\in \Z} P_{2i+1}}{}{}{\bigoplus_{i\in\Z
}P_{2i}}.\]
\end{lemma}

\begin{pf}
This is obvious (when you think about it).
\end{pf}

The functor $D$ is not an equivalence in general.  It is  an equivalence when $\A$ is hereditary \cite[Cor. 7.1]{peng_xiao}, although we shall not need this fact.

\subsection{Acyclic complexes}
\label{ac}
As stated above, from now on we consider only $\Zt$ graded complexes, and drop the $\Zt$ subscripts.
A complex $M_\blob\in\CC(\A)$ is called acyclic if $H_*(M_\blob)=0$.
For any object $P\in\P$, there are acyclic complexes
\[ K_P= \big(\cplx{P}{1}{0}{P}\big), \qquad K_P^\dag= \big(\cplx{P}{0}{1}{P}\big).\]
Note that a complex $M_\blob\in\CC(\P)$ is acyclic precisely if $M_\blob\isom 0$ in $\Ho(\P)$.

\begin{lemma}
\label{acyclic}
If  $M_\blob\in\CC(\P)$ is an acyclic complex of projectives, then
there are objects $P,Q\in \P$, unique up to isomorphism, such that
 $M_\blob\isom K_P\oplus K_Q^\dag$.
\end{lemma}

\begin{pf}
Suppose the complex
\[\cplx{M_1}{d_1}{d_0}{M_0}\]
is acyclic. Set $P=\im(d_0)=\ker(d_1)$ and $Q=\ker(d_0)=\im(d_1)$. There are short exact sequences
\begin{equation}
\label{l}0\lra P\lra M_1\lra Q\lra 0, \qquad 0\lra Q\lra M_0\lra P\lra 0.\end{equation}
The long exact sequence in $\Ext_\A(\blank,\blank)$ shows that for $i>0$
\[\Ext^i_\A(P,\blank)\isom \Ext^{i+1}_\A(Q,\blank)\isom \Ext^{i+2}_\A(P,\blank).\]
Since $\A$ has finite global dimension these groups must all vanish, and hence  $P$ and $Q$ are projective. The short exact sequences \eqref{l} can then be split, and it follows easily that $M_\blob\isom K_P\oplus K_Q^\dag$.
\end{pf}

\subsection{Extensions of complexes}
We shall make frequent use of the following simple result. 

\begin{lemma}
\label{exts}
If $M_\blob,N_\blob\in \CC(\P)$ then
\[\Ext^1_{\CC(\A)}(N_\blob,M_\blob)\isom \Hom_{\Ho(\A)}(N_\blob,M_\blob^\dag).\]
\end{lemma}

\begin{pf}
Consider a short exact sequence in $\CC(\A)$
\begin{equation}
\label{ter}
0\lra M_\blob\lra P_\blob\lra N_\blob\lra 0.\end{equation}
Since $M_\blob$ and $N_\blob$ are projective complexes, this is isomorphic as an extension in $\CC(\A)$ to a sequence of the form
\[\begin{xy}
\xymatrix{
M_1\ar[d]_{i_1}\ar@<.5ex>[r]^{f_M}		& \ar@<.5ex>[l]^{g_M} M_0\ar[d]^{i_0}\\
M_1\oplus N_1\ar[d]_{p_1}\ar@<.5ex>[r]^{f}	& \ar@<.5ex>[l]^{g}	 M_0\oplus N_0\ar[d]^{p_0}
\\
N_1\ar@<.5ex>[r]^{f_N}	& \ar@<.5ex>[l]^{g_N}	 N_0}
\end{xy} \]
where $i_0,i_1$ and $p_0,p_1$ denote the canonical inclusions and projections. Writing $f$ and $g$ in matrix form we have
\[f=\mat{f_M}{s_1}{0}{f_N}, \quad g=\mat{g_M}{s_0}{0}{g_N}\]
for morphisms $s_i\colon N_i\to M_{i+1}$. The condition that $f\circ g=g\circ f=0$ is precisely the condition that
\begin{equation}
\label{map}
s_\blob\colon N_\blob\to M_\blob^\dag\end{equation}
is a morphism of complexes.  Thus every morphism of complexes \eqref{map} determines an extension \eqref{ter}, and conversely, every extension \eqref{ter} arises in this way.

Suppose now that we have two morphisms \[s_\blob, \tilde{s}_\blob \colon N_\blob\to M_\blob^\dag\] determining complexes $P_\blob$ and $\tilde{P}_\blob$. The corresponding extension classes \eqref{ter} agree precisely if there is an  isomorphism $k_\blob\colon P_\blob\to \tilde{P}_\blob$ in $C(\A)$ commuting with the identity maps on $M_\blob$ and $N_\blob$. This last condition means that in matrix notation $k_\blob$ takes the special form
\[k_1=\mat{1}{h_1}{0}{1}, \quad k_0=\mat{1}{h_0}{0}{1}\]
for morphisms $h_i\colon N_i\to M_i$ in $\A$.
The condition that $k_\blob$ is a morphism of complexes is
precisely the condition that the diagram
\[\begin{xy}
\xymatrix{
M_1\oplus N_1\ar[d]_{k_1}\ar@<.5ex>[r]^{f}		& \ar@<.5ex>[l]^{g} M_0\oplus N_0\ar[d]^{k_0}\\
M_1\oplus N_1\ar@<.5ex>[r]^{\tilde{f}}	& \ar@<.5ex>[l]^{\tilde{g}}	 M_0\oplus N_0}
\end{xy} \]
commutes, which translates into the statement that
$h_\blob$
is a homotopy relating $s_\blob$ and $\tilde{s}_\blob$.
\end{pf}

\subsection{Hall algebra}
\label{burbage}

Let $\H(\CC(\A))$ be the Hall algebra of the abelian category $\CC(\A)$ as defined in Section \ref{hall}. Note that the definition makes perfect sense, even though $\CC(\A)$ is not usually of finite global dimension (even when $\A=\Vect_k$). Indeed, the spaces $\Ext^1_{\CC(\A)}(N_\blob,M_\blob)$ are all finite-dimensional by Lemma \ref{exts}.
Let \[\H(\CC(\P))\subset \H(\CC(\A))\] be the subspace spanned by complexes of projective modules.
It is a subalgebra under the Hall product because the subcategory $\P(\A)\subset\CC(\A)$ is closed under extensions. 

Define  $\Htw(\CC(\P))$ to be the same vector space as $\H(\CC(\P))$, but with twisted  multiplication
\begin{equation}
\label{tw}\gen{M_\blob}*\gen{N_\blob}=t^{\<M_0,N_0\>+\<M_1,N_1\>}\cdot \gen{M_\blob}\diam \gen{N_\blob}.\end{equation}
This is a slight abuse of terminology,  since the multiplication  is not being twisted  by the Euler form of the category $\CC(\A)$. Indeed,  since $\CC(\A)$ usually has infinite global dimension,  it does not have a well-defined Euler form.

\smallskip

The  acylic complexes $K_P$ that were introduced in Section \ref{ac} define elements of $\Htw(\CC(\P))$ with particularly simple properties.

\begin{lemma}
\label{easy}
For any object $P\in \P$ and any complex $M_\blob\in\CC(\P)$ there are identities in $\Htw(\CC(\P))$
\begin{equation}\gen{K_P}*\gen{M_\blob}= t^{\<\cl{P},\cl{M_\blob}\>}\cdot\gen{K_P\oplus M_\blob},\end{equation}\begin{equation}\gen{M_\blob}*\gen{K_P}= t^{-\<\cl{M_\blob},\cl{P}\>}\cdot \gen{ K_P\oplus M_\blob}.\end{equation}
\end{lemma}

\begin{pf}
It is easy to check directly that
\[\Hom_{\CC(\A)}(K_P,M_\blob)=\Hom_\A(P,M_1), \quad \Hom_{\CC(\A)}(M_\blob,K_P)=\Hom_\A(M_0,P).\]
The complexes $K_P$ are homotopy equivalent to the zero complex, so Lemma \ref{exts} shows that the extension group in the definition of the Hall product vanishes. Taking into account the twisting \eqref{tw} gives the result.
\end{pf}

\begin{lemma}
\label{cor}
For any object $P\in \P$ and any complex $M_\blob\in\CC(\P)$ there are identities in $\Htw(\CC(\P))$
\begin{equation}\label{eq1} \gen{K_P}*\gen{M_\blob}=t^{(\cl{P},\cl{M_\blob})} \cdot\gen{M_\blob}* \gen{K_P}\end{equation}\begin{equation}\label{eq2}
\gen{K_P^\dag}*\gen{M_\blob}=t^{-(\cl{P},\cl {M_\blob})} \gen{M_\blob}*\cdot\gen{K_P^\dag}.\end{equation}
\end{lemma}

\begin{pf}
Equation \eqref{eq1} is immediate from Lemma \ref{easy}. Equation \eqref{eq2} follows  by applying the involution $\dag$.
\end{pf}

In particular, since $\cl{K_{P}}=0\in K(\A)$, we have
\begin{equation}
\label{additive}
\gen{K_P} * \gen{K_Q} = \gen{K_{P}\oplus K_{Q}}, \quad \gen{K_P}*\gen{K_Q^\dag}=\gen{K_P\oplus K_Q^\dag}.\end{equation}
\begin{equation}
[\gen{K_P},\gen{K_Q}]=[\gen{K_P},\gen{K_Q^\dag}]=[\gen{K_P^\dag},\gen{K_Q^\dag}]=0. \end{equation}
Note that  any element of the form $\gen{K_{P}}* \gen{K_P^\dag}$ is central.
%

\subsection{Localization and reduction}
\label{def}

%
%
%
%
%

Define the localized Hall algebra $\DH(\A)$ to be the localization of $\Htw(\CC(\P))$ with respect to the elements $\gen{M_\blob}$ corresponding to acyclic complexes $M_\blob$. In symbols
\[\DH(\A)=\Htw(\CC(\P))[ \gen{M_\blob}^{-1} : H_*(M_\blob)=0].\]
By \eqref{additive} and Lemma \ref{acyclic} this is the same as localizing by the elements $\gen{K_P}$ and $\gen{K_P^\dag}$ for all objects $P\in \P$. 
Note that Lemma \ref{cor} shows that these elements satisfy the Ore conditions, and hence give a well-defined localization. 

The assignment $P\mapsto K_P$ extends to a  group homomorphism
\[K\colon K(\A)\lra \DH(\A)^\times.\]
Explicitly this is given by writing an element $\alpha\in K(\A)$ in the form  $\alpha=\cl{P}-\cl{Q}$ for objects $P,Q\in \P,$ and  setting
\[K_\alpha=K_{\cl{P}}* K_{\cl{Q}}^{-1}.\]
Composing with the involution $\dag$ gives another map
\[K^\dag\colon K(\A)\lra \DH(\A)^\times.\]
Note that \eqref{eq1} and \eqref{eq2} continue to hold with the elements $\gen{K_P}$ and $\gen{K_P^\dag}$ replaced by $K_\alpha$ and $K_\alpha^\dag$  for arbitrary $\alpha\in K(\A)$.

 Define the reduced localized Hall algebra by setting $\gen{M_\blob}=1$ whenever $M_\blob$ is an acyclic complex, invariant under the shift functor.  In symbols
\[\DHred(\A)=\DH(\A)/\big(\gen{M_\blob}-1: H_*(M_\blob)=0, \  M_\blob\isom M_\blob^\dag\big).\]  
By Lemma \ref{acyclic} this is the same as setting
\begin{equation*}
\label{hot}\gen{K_{P}}* \gen{K_P^\dag} = 1\end{equation*}
for all $P\in \P$.  
Since these  relations force $K_P$ to be invertible, we also have
\[\DH_{\red}(\A)=\Htw(\P))/\big(\gen{M_\blob}-1: H_*(M_\blob)=0, \  M_\blob\isom M_\blob^\dag\big).\]
We leave it to the reader to check that the shift functor $\dag$ defines involutions of both $\DH(\A)$ and $\DH_{\red}(\A)$.

%
%
%
%

\section{Hereditary case}
\label{her}

In this section we prove our main theorem realizing quantum enveloping algebras as Hall algebras of complexes. Throughout $\A$ is an abelian category satisfying the assumptions (a)--(e) of Section \ref{ass}. In particular $\A$ is hereditary, which means that every subobject of a projective object is also projective.

\subsection{Minimal resolutions}
The conditions on $\A$ imply that every object $A\in \A$ has a projective resolution of the form
\begin{equation}
\label{res}
0\lra P\lRa{f} Q \lRa{g} A\lra 0.\end{equation}
 Condition (a) ensures that the category $\A$ has the Krull-Schmidt property\footnote{Reference required}. Decomposing $P$ and $Q$ into finite direct sums
$P=\oplus P_i$, $Q=\oplus Q_j$,
we may write $f=(f_{ij})$ in matrix form for certain morphisms $f_{ij}\colon P_i\to Q_j$. The resolution \eqref{res} is said to be minimal if none of the morphisms $f_{ij}$ is an isomorphism. It is easy to see that any object has a minimal resolution. In fact,

\begin{lemma}
\label{fag}
Any resolution \eqref{res} is isomorphic to a resolution of the form
\[0\lra R\oplus P\lRa{1\oplus f'} R\oplus Q'\lRa{(0,g')}A\lra 0,\]
for some object $R\in \P$, and some minimal projective resolution
\[0\lra P'\lRa{f'} Q' \lRa{g'} A\lra 0.\]
\end{lemma}

\begin{pf}
If \eqref{res} is not minimal, some $f_{ij}$ is an isomorphism. Without loss of generality we can assume that $P_i=Q_j$ and $f_{ij}$ is the identity. Set $P'=P/P_i$ and $Q'=Q/Q_j$. Then there is a split short exact sequence of complexes
\[\begin{xy} \xymatrix@=0pt@C=20pt@R=15pt{
R \ar[d]\ar[r]^{1}& R\ar[d] \\
P \ar[r]^{f}\ar[d] &Q \ar[r]^{g}\ar[d] &A\ar[d]^{1} \\
P'\ar[r]^{f'} &Q'\ar[r]^{g'} &A\\}\end{xy}\]
If we choose the splittings correctly we get $f=1\oplus f'$. We can now repeat. This process eventually terminates by the finiteness assumption.
\end{pf} 
 
 \subsection{The complexes $C_A$}
Given an object $A\in \A$, take a minimal projective resolution \eqref{res} and
consider the corresponding $\Zt$-graded complex
\begin{equation}
\label{cp}C_A=(\cplx{P}{f}{0}{Q})\in \CC(\P).\end{equation}
Lemma \ref{fag}  shows that any two minimal projective resolutions of  $A$ are isomorphic, so the complex $C_A$ is well-defined up to isomorphism.  The following result seems to be well known, but I couldn't find a proof in the literature.

\begin{lemma}
\label{decomp}
Every object $M_\blob\in C(\P)$ has a direct sum decomposition
\[ M_\blob = C_A \oplus C_B^\dag\oplus K_P\oplus K_Q^\dag.\]
Moreover, the objects $A,B\in \A$ and $P,Q\in \P$ are uniquely determined up to isomorphism.
\end{lemma}

\begin{pf}
By the Krull-Schmidt property we can assume that our complex
\[\cplx{M_1}{d_1}{d_0}{M_0}\]
 is indecomposable. Consider the short exact sequences
\[0\lra \ker(d_1)\lRa{i} M_1\lRa{p} \im(d_1)\lra 0,\]\[ 0\lra \ker(d_0)\lRa{j} M_0\lRa{q} \im(d_0)\lra 0.\]
By the hereditary assumption, all the objects appearing in these sequences are projective. Thus the sequences split, and we can find morphisms
\[r\colon M_1\to \ker(d_1), \quad k\colon \im(d_0)\to M_0\]
such that $r\circ i=\id$ and $q\circ k=\id$. Now there are morphisms of complexes
\[\begin{xy}
\xymatrix{
M_1\ar[d]_{r}\ar@<.5ex>[r]^{d_1}	& \ar@<.5ex>[l]^{d_0}	 M_0\ar[d]^{q}
\\
\ker(d_1)\ar@<.5ex>[r]^{0}	& \ar@<.5ex>[l]^{m}	 \im(d_0)}
\end{xy} \qquad \begin{xy}
\xymatrix{
M_1\ar@<.5ex>[r]^{d_1}	& \ar@<.5ex>[l]^{d_0}	 M_0
\\
\ker(d_1)\ar[u]^{i}\ar@<.5ex>[r]^{0}	& \ar@<.5ex>[l]^{m}	 \im(d_0)\ar[u]_{k}}
\end{xy}\]
where $m$ is the obvious inclusion (note that $d_0=i\circ m\circ q$). Since $M_\blob$ is indecomposable we can conclude that either \[\ker(d_1)=\im(d_0)=0\quad \text{or}\quad \ker(d_0)=\im(d_1)=0.\] Without loss of generality we assume the first holds, so that $d_0=0$ and
 $d_1$ is injective. Then
\[0\lra P\lRa{d_1} Q\lra A\lra 0\]
is a projective resolution of $A=H_0(M_\blob)$. Since  $M_\blob$ is assumed indecomposable, Lemma \ref{fag} shows that  either $M_\blob\isom K_P$ or   $M_\blob\isom C_{A}$.
\end{pf}

\subsection{Embedding $\H(\A)$ in $\DH(\A)$}

Given an object $A\in \A$ define 
\[\E_A=  t^{\<\cl{P},\cl{A}\>}\cdot K_{-\cl{P}}*\gen{ C_A}\in \DH(\P).\]
An explanation for the strange-looking coefficients is as follows. Suppose  we take a different, not necessarily minimal, projective resolution \eqref{res}, and consider the corresponding complex \eqref{cp}  in $\CC(\P)$.  By Lemma \ref{fag} this will have the effect of replacing the complex $C_A$ by a complex  $K_R\oplus C_A$ for some object $R\in \P$. But then, by Lemma \ref{easy}
\[ t^{\<\cl{P\oplus R},\cl{A}\>}\cdot K_{-\cl{P\oplus R}}*\gen{  K_R\oplus C_A}
=t^{\<\cl{P},\cl{A}\>}\cdot K_{-\cl{P}}*\gen{ C_A},\]
so we get the same element $E_A$.

\begin{lemma}
\label{embed}
There is an injective ring homomorphism
\[I_+\colon \Htw(\A)\longinto \DH(\A), \quad [A] \mapsto \E_A.\] 
\end{lemma}

\begin{pf}
Given objects $A_1,A_2\in \A$, take minimal projective resolutions
\[0\lra P_i\stackrel{f_i}{\lra}Q_i\lra A_i\lra 0,\]
and define objects $C_{A_i}\in \CC(\P)$ as before.
First note that, by Lemma \ref{easy}  and Lemma \ref{cor},
\[\E_{A_1}* \E_{A_2}=t^{\<\cl{P_1},\cl{A_1}\>+\<\cl{P_2},\cl{A_2}\>}\cdot K_{-\cl{P_1}}*\gen{ C_{A_1}}*K_{-\cl{P_2}}*\gen{ C_{A_2}}\]
\[=t^{(\cl{P_2},\cl{A_1})+\<\cl{P_1},\cl{A_1}\>+\<\cl{P_2},\cl{A_2}\>}\cdot K_{-(\cl{P_1}+\cl{P_2})}*\gen{ C_{A_1}}*\gen{ C_{A_2}}.\]
Since the complexes $C_{A_i}$ are quasi-isomorphic to the objects $A_i$, Lemmas \ref{missing} and \ref{exts} implies that
\[\Ext^1_{\CC(\A)}(C_{A_1},C_{A_2})= \Hom_{\Ho(\A)}(C_{A_1},C_{A_2}^\dag) =\Ext^1_\A(A_1,A_2).\]
Furthermore, any extension of $C_{A_1}$ by $C_{A_2}$ is the complex  $C_{A_3}$ defined by the corresponding extension $A_3$ of $A_1$ by $A_2$. 
To see this, consider the diagram
\[\begin{xy}
\xymatrix@=0pt@C=20pt@R=15pt{
P_2\ar[d]_{i}\ar@<.5ex>[r]^{f_2}		& \ar@<.5ex>[l]^{0} Q_2\ar[d]^{i}\\
P_1\oplus P_2\ar[d]_{p}\ar@<.5ex>[r]^{u}	& \ar@<.5ex>[l]^{v}	 Q_1\oplus Q_2\ar[d]^{p}
\\
P_1\ar@<.5ex>[r]^{f_1}	& \ar@<.5ex>[l]^{0}	 Q_1}
\end{xy} \]
Since the $f_i$ are monomorphisms so is $u$. But then  $u\circ v=0$ implies that $v=0$.

It is easy to check directly that there is a short exact sequence
\[0\lra \Hom_\A(Q_1,P_2) \lra \Hom_{\CC(\A)}(C_{A_1},C_{A_2})\lra \Hom_\A(A_1,A_2)\lra 0.\]
Putting all this together  we get
\[\E_{A_1}* \E_{A_2}=t^n\cdot \big.\sum_{A_3\in \Iso(\A)} \frac{|\Ext^1_\A(A_1,A_2)_{A_3}|}{|\Hom_\A(A_1,A_2)|}\cdot E_{A_3},\]
where, taking also the twisting into account, the total power of $t$  is
\[n=\<\cl{P_2},\cl{A_1}\>+\<\cl{A_1},\cl{P_2}\>+\<\cl{P_1},\cl{A_1}\>+\<\cl{P_2},\cl{A_2}\>-\<\cl{P_1}+\cl{P_2}, \cl{A_1}+\cl{A_2}\>\]\[+ \<\cl{P_1},\cl{P_2}\>+\<\cl{Q_1},\cl{Q_2}\>-2\<\cl{Q_1},\cl{P_2}\>.\]
Using $\cl{Q_i}=\cl{P_i}+\cl{A_i}$ this reduces to 
$\<\cl{A_1},\cl{A_2}\>$
as required.

To show that the map $I_+$ is injective note that the linear map
\[Q\colon \DH(\A)\to \Htw(\A), \quad \gen{M_\blob}\mapsto t^{-\<\cl{M_1}, \cl{H_0(M_\blob)}\>}\cdot \gen{H_0(M_\blob)}.\]
satisfies $Q( I_+ (\gen{A}))= \gen{A}$.
\end{pf}

\subsection{Commutation relations}
Composing the map $I_+$ of Lemma \ref{embed} with the involution $\dag$ gives another injective ring homomorphism
\[I_-\colon \Htw(\A)\to \DH(\A), \quad [A] \mapsto F_A,\]
where for any object $A\in \A$ we set $F_A=E_A^\dag$.

\begin{lemma}
\label{commute}
Suppose $A_1,A_2\in \A$ satisfy
\begin{equation}\label{assump}\Hom_\A(A_1,A_2)=0=\Hom_\A(A_2,A_1).\end{equation}
Then
$[E_{A_1} ,F_{A_2}]=0$.
\end{lemma}

\begin{pf}
Take minimal projective resolutions
\[0\lra P_i\stackrel{f_i}{\lra}Q_i\lra A_i\lra 0\] as before.
Since the complexes $C_{A_i}$ are quasi-isomorphic to the objects $A_i$ it follows from Lemmas \ref{missing} and \ref{exts} that
\[\Ext^1_{\CC(\A)}(C_{A_1},C_{A_2}^\dag)=\Hom_{\Ho(\A)}(C_{A_1},C_{A_2})=\Hom_\A(A_1,A_2).\]
It is easy to check directly that
\begin{equation}
\label{eyes}\Hom_{\CC(\A)}(C_{A_1},C_{A_2}^\dag)=\Hom_\A(P_1,Q_2).\end{equation}
Thus we have
\[\gen{C_{A_1}}*\gen{C_{A_2}^\dag}= t^{\<\cl{P_1},\cl{Q_2}\>+\<\cl{Q_1},\cl{P_2}\>}\cdot q^{-\<\cl{P_1},\cl{Q_2}\>} \cdot \gen{C_{A_1}\oplus C_{A_2}^\dag},\] 
and hence
\begin{gather*}\E_{A_1}* F_{A_2}=
t^{\<\cl{P_1},\cl{A_1}\>+\<\cl{P_2},\cl{A_2}\>}\cdot K_{-\cl{P_1}}*\gen{ C_{A_1}}*K_{-\cl{P_2}}^\dag*\gen{ C_{A_2}^\dag}\\
= t^n\cdot K_{-\cl{P_1}} * K_{-\cl{P_2}} ^\dag* \gen{C_{A_1}\oplus C_{A_2}^\dag},\end{gather*}
where the total power of $t$ is
\[n=\<\cl{Q_1},\cl{P_2}\>-\<\cl{P_1},\cl{Q_2}\>+\<\cl{P_1},\cl{A_1}\>+\<\cl{P_2},\cl{A_2}\>-\<\cl{A_1},\cl{P_2}\>-\<\cl{P_2},\cl{A_1}\>.\]
Using $\cl{Q_i}=\cl{P_i}+\cl{A_i}$ this can be rewritten
\begin{equation}\label{n} n=\<\cl{P_1},\cl{A_1}\>+\<\cl{P_2},\cl{A_2}\>-\<\cl{P_1},\cl{A_2}\>-\<\cl{P_2},\cl{A_1}\>,\end{equation}
which is invariant under exchanging 1 and 2.
Then
\begin{gather*}F_{A_2}*E_{A_1}=(E_{A_2}*F_{A_1})^\dag= t^{ n} \cdot K_{-\cl{P_2}}^\dag* K_{-\cl{P_1}} * \gen{C_{A_2}\oplus C_{A_1}^\dag}^\dag\\
= t^{ n} \cdot K_{-\cl{P_1}} * K_{-\cl{P_2}}^\dag * \gen{C_{A_1}\oplus C_{A_2}^\dag}=E_{A_1}* F_{A_2}\end{gather*}
as claimed.
\end{pf}

To complete the derivation of the quantum group relations we need one more computation.

\begin{lemma}
\label{tricky}
Suppose $A\in\A$ satisfies $\operatorname{End}_\A(A)=k$.
Then
\[[E_A ,F_A] =  (q-1)\cdot (K_{\cl{A}}^\dag - K_{\cl{A}}).\]
\end{lemma}

\begin{pf}
Take a projective resolution \eqref{res} as before. By Lemmas \ref{missing} and  \ref{exts} we have
\[\Ext^1_{\CC(\A)}(C_A,C_A^\dag)=\Hom_\A(A,A)=k.\]
Note that any endomorphism of $A$ is either an isomorphism or zero. It follows from the long exact sequence in cohomology that if
\[0\lra C_A^\dag \lra M_\blob \lra C_A\lra 0,\]
is a non-split short exact sequence in $\CC(\A)$ then $M_\blob$ is acyclic. Moreover it is easy to check that the kernel of the differential of $M_\blob$ must be isomorphic to $P$ in degree 0 and $Q$ in degree 1. Lemma \ref{acyclic} then implies that $M_\blob\isom K_P\oplus K_Q^\dag$.
Thus, by \eqref{eyes}
\[\gen{C_A}*\gen{C_A^\dag}= t^{\<\cl{Q},\cl{P}\>-\<\cl{P},\cl{Q}\>} \cdot \big( \gen{C_A\oplus  C_A^\dag} + (q-1) \cdot \gen{K_{P} \oplus K_{Q}^\dag}\big) .\]
Using the calculation given in the proof of Lemma \ref{commute}, and noting that $n=0$ when $A_1=A_2$, we have
\[E_A * F_A =   K_{-\cl{P}}*K_{-\cl{P}}^\dag *\big(\gen{C_A\oplus  C_A^\dag} + (q-1)\cdot  K_{\cl{P}}* K_{\cl{Q}}^\dag\big).\]
Applying $\dag$, subtracting and using $\cl{A}=\cl{Q}-\cl{P}$ gives the result.
\end{pf}

\subsection{Decompositions}
In this section we make some more precise statements about the relationships between the various Hall algebras we have been considering. 

\begin{lemma}
\label{embede}
 There is an embedding of algebras \[I_+^{\e}\colon \H^{\e}_{\tw}(\A)\longinto \DH(\A)\]
defined on generators by $\gen{A}\mapsto E_A$ and $K_{\alpha}\mapsto K_\alpha$.
\end{lemma}

\begin{pf}
The existence of such a ring homomorphism follows from Lemma \ref{embed} together with a comparison of the relations \eqref{twisted} defining the extended Hall algebra, with the relation \eqref{eq1}.
\end{pf}

Composing with the involution $\dag$ gives another such embedding
\[I_-^{\e}\colon\H^{\e}_{\tw}(\A)\longinto \DH(\A)\]
defined on generators by $\gen{B}\mapsto F_B$ and $K_{\beta}\mapsto K_\beta^\dag$.

\begin{lemma}
\label{ident}
The multiplication map $\mu\colon a\tensor b \mapsto I_+^{\e}(a) * I_-^{\e}(b)$ defines  an isomorphism of vector spaces 
\[\mu\colon \H^{\e}_{\tw}(\A) \tensor_\C \H^{\e}_{\tw}(\A) \lra\DH(\A).\]
\end{lemma}

\begin{pf}
Using Lemmas \ref{easy} and \ref{decomp} it is easy to see that the algebra $\DH(\A)$ has a basis consisting of elements \[\gen{C_A\oplus C_B^\dag}* K_\alpha * K^\dag_\beta, \quad A,B\in \Iso(\A), \quad \alpha,\beta\in K(\A).\]
Recall the partial order on $K(\A)$ defined in Section \ref{ass}, and define subspaces $\DH_{\leq \gamma}(\A)$ spanned by elements from this basis for which $\cl{A}+\cl{B}\leq \gamma$. We claim that
\begin{equation}
\label{pesto}
\DH_{\leq \alpha}(\A)*\DH_{\leq \beta}(\A)\subseteq \DH_{\leq \alpha+\beta}(\A)\end{equation}
so that this defines a filtration on $\DH(\A)$.

Consider an extension in $\CC(\P)$ 
\[0\lra M_\blob\lra P_\blob\lra N_\blob\lra 0.\]
The long exact sequence in homology can be split to give two long exact sequences
\begin{gather*}
0\lra K\lra H_0(M_\blob)\lra H_0(P_\blob)\lra H_0(N_\blob)\lra Q\lra 0, \\
0\lra Q\lra H_1(M_\blob)\lra H_1(P_\blob)\lra H_1(N_\blob)\lra K\lra 0.
\end{gather*}
It follows  that there is a relation in $K(\A)$
\[\cl{H_0(P_\blob)}+\cl{H_1(P_\blob)}+2(\cl{K}+\cl{Q})
=\cl{H_0(M_\blob)}+\cl{H_1(M_\blob)}+\cl{H_0(N_\blob)}+\cl{H_1(N_\blob)}\]
which proves \eqref{pesto}.

Suppose now that $N_\blob=C_A$ and $M_\blob=C_B^\dag$ for given objects $A,B\in \A$. Then $K=0$, and by Lemmas \ref{missing} and \ref{exts}
\[\Ext^1_{\CC(\A)}(N_\blob,M_\blob)=\Hom_\A(A,B),\]
and the extension class is completely determined by the connecting morphism
$H_0(N_\blob)\to H_1(M_\blob)$.  By assumption (e) of Section \ref{ass}, we therefore know that $\cl{Q}=0$ exactly when  the extension is trivial.
It follows that in the  graded algebra associated to the filtered algebra $\DH(\A)$, one has
a relation
\[\mu(\gen{A} *K_\alpha \tensor \gen{B}\tensor K_\beta)=t^n \cdot \gen{C_{A}\oplus C_B^\dag} *K_{\alpha-\cl{N}_1} * K_{\beta-\cl{M}_0}^\dag,\]
where $n$ is some integer.
It follows that $\mu$ takes a basis to a basis and is hence  an isomorphism. 
\end{pf}

We shall also need the following triangular decomposition statement.

\begin{lemma}
\label{tri2}
The multiplication map $\gen{A}\tensor \alpha\tensor \gen{B}\mapsto E_A * K_\alpha * F_B$ defines
 an isomorphism of vector spaces
\[H_{\tw}(\A)\tensor_{\C} \C[K(\A)]\tensor_\C \H_{\tw}(\A) \lra \DH_{\red}(\A).\]
\end{lemma}

\begin{pf}
The same argument given for Lemma \ref{ident} also applies here.
\end{pf}

\subsection{Quantum enveloping algebras}

\label{quantum}

 Let $Q$ be a finite quiver without oriented cycles, and vertices $\{1,\cdots, n\}$. Let $\A$ be the category of finite-dimensional representations of $Q$ over the field $k=\F_q$. This abelian category satisfies all the assumptions (a)--(e) of Section \ref{ass}.
For each vertex $i$ of $Q$ we denote the corresponding one-dimensional simple module by $S_i\in\A$.
Note that \[a_{ij}=(S_i,S_j)\] is the Cartan matrix of the graph underlying $Q$.

\begin{thm}There is an injective homomorphism of algebras
\[R\colon \U_t(\g)\longinto\DH_{\red}(\A),\]
defined on generators by
\[R(E_i)=(q-1)^{-1}\cdot E_{S_i}, \quad R(F_i) = (-t)\cdot (q-1)^{-1}\cdot {F_{S_i}},\]
\[R(K_i)= K_{\cl{S_i}}, \quad R(K_i^{-1}) =K_{\cl{S_i}}^\dag.\]
The map $R$ is an isomorphism precisely if the  graph underlying $Q$ is a simply-laced Dynkin diagram.
\end{thm}

\begin{pf}
The existence of the ring homomorphism follows from the results of the previous sections, which show that the given elements satisfy the defining relations of the quantum group. There is a commutative diagram of linear maps
\[\begin{xy}
\xymatrix{ \U_t(\n^+)\tensor\U_t(\h)\tensor \U_t(\n^-)\ar[rr]^{A}\ar[d]&&\H_{\tw}(\A)\tensor\C[K(\A)]\tensor \H_{\tw}(\A) \ar[d]\\
\U_t(\g) \ar[rr]^{R} &&\DH_{\red}(\A)}
\end{xy} \]
in which the vertical arrows are the isomorphisms described in Lemmas \ref{tri} and \ref{tri2}, and the homomorphism $A$ is built out of the homomorphisms of Theorem \ref{ring} in the obvious way. It follows that $R$ is injective in general, and an isomorphism precisely in the Dynkin case.
\end{pf}

Note that the definition of the map $R$ is slightly asymmetric. If we omit the factor $(-t)$ in the definition of $R(F_i)$ we get an extra factor of $(-t^{-1})$ on the right hand side of equation  \eqref{anom}. From the point of view of Hall algebras it would be more natural to write the relations in this way. Compare \cite[Remark 3.12]{burban}.

\bibliographystyle{amsplain}

\end{document}